# Remarks regarding the gap between continuous, Lipschitz, and differentiable storage functions for dissipation inequalities appearing in $H_\infty$ control

Lionel Rosier
Laboratoire d'Analyse Numérique et EDP,
Université Paris 11, bât. 425,
F91405 Orsay Cedex, France
Lionel.Rosier@math.u-psud.fr
and
Eduardo D. Sontag *
Department of Mathematics
Rutgers University
New Brunswick, NJ 08903, USA
sontag@control.rutgers.edu

**Abstract**

This paper deals with the regularity of solutions of the Hamilton-Jacobi Inequality which arises in $H_\infty$ control. It shows by explicit counterexamples that there are gaps between existence of continuous and locally Lipschitz (positive definite and proper) solutions, and between Lipschitz and continuously differentiable ones. On the other hand, it is shown that it is always possible to smooth-out solutions, provided that an infinitesimal increase in gain is allowed.

## 1 Introduction

The so-called "$H_\infty$ control problem" is that of finding a state (or more generally, a measurement-based) feedback which stabilizes a given system, while satisfying an energy-gain ($L^2$ operator norm) constraint. For linear systems, the problem has a long history, and an elegant solution was provided in the by now classical paper [8]. The nonlinear version of this problem has been the subject of intense research as well; see for instance [3, 4, 13, 20] among many others. A central role in these studies is played by a partial differential inequality, a Hamilton-Jacobi Inequality (HJI) which is satisfied by a "storage" or "energy" function $V$ associated with the closed loop system.

In this paper, we concern ourselves with the analysis of the HJI for a system already in closed-loop form, since we wish to make some remarks about the regularity of solutions of this equation. Moreover, in order to keep the discussion as simple as possible, we analyze the case of full state measurements, but similar conclusions could be drawn for the case when outputs are of interest. Thus, the main focus of this paper will be on systems, affine in inputs, of the following form:

$$\dot{x} \;=\; g_0(x) \;+\; \sum_{i=1}^{m} u_i\, g_i(x) \tag{1}$$

for which states $x$ evolve in $\mathbb{R}^n$, and inputs $u = (u_1, \ldots, u_m)$ take values in $\mathbb{R}^m$. We assume that the vector fields $g_i$, $i = 0, \ldots, m$ are locally Lipschitz in $x$. (These systems might be thought of having been obtained from a more general class of systems $\dot{x} = f(x, u, v)$ after applying a stabilizing feedback $v = k(x)$, but that interpretation is irrelevant to the results that we give.) In this introduction, we restrict attention to the above class of systems, but later we will also provide several results valid for more general systems, not necessarily affine in inputs, of the form

$$\dot{x} \;=\; F(x, u)\,, \tag{2}$$

with $x(t) \in \mathbb{R}^n$ and $u(t) \in \mathcal{U}$. (As far as the definitions are concerned, we do not need impose any technical conditions on $\mathcal{U}$ and $F$.) We use $|\cdot|$ to denote Euclidean norm in $\mathbb{R}^n$ or $\mathbb{R}^m$. For systems (1), the HJI's of interest are inequalities which are often expressed as:

$$\nabla V(x) \cdot g_0(x) \;+\; \frac{1}{4\gamma} \sum_{i=1}^{m} (\nabla V(x) \cdot g_i(x))^2 \;+\; |x|^2 \;\leq\; 0\,, \tag{3}$$

where $\gamma > 0$ is the (square of the) "$L^2$ gain" of the system. Inequality (3) for non-differentiable $V$ must be interpreted in a generalized sense, as we discuss below. The $V$'s that arise, because of implicit stability considerations, must be proper and positive definite.

---

*This research was supported in part by Grant F49620-98-1-0242.

Modulo some elementary calculus to do the implied maximization over $u$, the inequality (3) is equivalent to:

$$\nabla V(x) \cdot \left( g_0(x) \;+\; \sum_{i=1}^{m} u_i\, g_i(x) \right) \;\leq\; -|x|^2 + \gamma |u|^2\,, \tag{4}$$

understood as holding *for all u*. This latter form is preferred, because it makes sense for arbitrary (not necessarily input-affine) systems (2), and because the theory then makes a natural contact with the theory of dissipative systems developed by Willems and Moylan and Hill (see e.g. [11, 12, 20, 21]). Under suitable technical assumptions, one can show that (4) is equivalent to the following inequality:

$$V(x(b)) - V(x(a)) \;\leq\; \int_a^b |u(t)|^2 - |x(t)|^2\, dt \tag{5}$$

holding along all solutions $(x(\cdot), u(\cdot))$ of (1) with $x(\cdot)$ absolutely continuous and $u(\cdot)$ measurable and essentially bounded, for each $a < b$ in the domain of the solution. This follows from general results about proximal characterizations of "strong invariance" of differential inclusions given in Chapter 4 of [6], or the equivalent viscosity characterizations given e.g. in [14], see also the discussion in [3]. Notice that (5) has the following consequence for trajectories which start at $x(0) = 0$, and if $V$ is nonnegative: for all $T > 0$ for which the solution is defined, it must hold that

$$\int_0^T |x(t)|^2\, dt \;\leq\; \gamma \int_a^b \gamma |u(t)|^2\, dt\,.$$

This means that the map $u(\cdot) \mapsto x(\cdot)$ seen as a map between square-integrable functions, has operator norm $\leq \sqrt{\gamma}$. Conversely, the theory relates such operator norms back to the existence of $V$'s solving the differential inequality. We will not discuss this well-known material any further, but rather consider the differential inequality as our object of study.

In general, it is not natural to impose the requirement that solutions to (4) should be smooth, so the HJI must be interpreted in the viscosity sense (a "viscosity supersolution"), see e.g. [7, 9]), or, in what is an essentially equivalent manner, the proximal analysis formalism found in [6]. Under suitable controllability hypotheses, it does make sense to restrict attention to *continuous* $V$'s, see for instance [3]. Thus we will always assume in our study that $V$ is at least continuous.

The main results in this note address the gap between continuous and $C^1$ solvability. We give examples which show (a) that there may exist even Lipschitz-continuous (proper and positive definite) solutions, but no possible continuously differentiable ones, and (b) there may exist continuous (proper and positive definite) solutions but no possible locally Lipschitz continuous ones. On the other hand, we provide a smoothing result which shows that it is possible to pass from $C^0$ to $C^1$ solutions as long as an infinitesimal increase in the gain $\gamma$ is allowed. A last section treats the special case of one-dimensional systems, where no gap exists. None of these facts is unexpected, of course, but it would appear that the corresponding counterexamples and proofs are not available in the literature.

## 2 Definitions and Statements of Results

Recall (cf. [6], Section 3.4) that a vector $\zeta \in \mathbb{R}^n$ is a *viscosity subgradient* of a function $V : \mathbb{R}^n \to \mathbb{R}$, at the point $x \in \mathbb{R}^n$, if

$$\liminf_{h \to 0} \frac{1}{|h|} \left[ V(x+h) - V(x) - \zeta \cdot h \right] \geq 0 \qquad (6)$$

The (possibly empty) set of all viscosity subgradients of $V$ at $x$ is called the *viscosity subdifferential*, and is denoted $\partial_D V(x)$. Observe that, if the function $V$ is differentiable at $x$, then $\partial_D V(x) = \{\nabla V(x)\}$.

**Definition 2.1** Suppose given a system $\Sigma$ as in (2), and a $\gamma \geq 0$. We say that a function $V : \mathbb{R}^n \to \mathbb{R}_{\geq 0}$ *witnesses the gain* $\gamma$ if the following condition holds:

$$\zeta \cdot F(x,u) \leq -|x|^2 + \gamma |u|^2 \quad \forall \zeta \in \partial_D V(x) \qquad (7)$$

for all $x \in \mathbb{R}^n \setminus \{0\}$ and $u \in \mathcal{U}$. The set of all those continuous $V$ which witness a given gain $\gamma$ for the system $\Sigma$ is denoted as $\mathcal{W}(\Sigma, \gamma)$. □

We recall that a continuous function $V : \mathbb{R}^n \to \mathbb{R}_{\geq 0}$ is said to be *proper* (or "radially unbounded") provided that every set of the form $\{x \,|\, V(x) \leq a\}$ is compact, for each $a > 0$ (equivalently, $V(x) \to \infty$ as $|x| \to \infty$), and is said to be *positive definite* if $V(x) = 0$ if and only if $x = 0$. We denote

$$\mathcal{W}_\infty(\Sigma, \gamma) := \{V \in \mathcal{W}(\Sigma, \gamma) \,|\, V \text{ is proper and positive definite}\}.$$

When $V$ is $\mathcal{C}^1$, condition (7) means simply that $\nabla V(x) \cdot F(x,u) \leq -|x|^2 + \gamma |u|^2$ for all $x \neq 0$ and $u$. As is well-known, asking for a globally $\mathcal{C}^1$ such function which is also positive

definite is overly restrictive. To see this, consider an input-affine system with $n = 1$, and take $u = 0$. The inequality would force the bound $|x|^2/|g_0(x)| \leq |V'(x)|$ for all $x \neq 0$. If $V$ is positive definite, it must have a local minimum at zero, so $V'(0) = 0$. Hence

$$|x|^2 = o(g_0(x)) \text{ as } x \to 0,$$

which is too strong a constraint on $g_0$ (one would not be able to study a system such as $\dot{x} = -x^3$). Therefore, it is routine (see for instance [19]) to drop the requirement of differentiability at the origin. We will denote by $\mathcal{C}_0^1$ the set of continuous functions $V : \mathbb{R}^n \to \mathbb{R}_{\geq 0}$ whose restriction to $\mathbb{R}^n \setminus \{0\}$ is continuously differentiable.

The first result consists of an example which shows that there may exist a globally Lipschitz $V \in \mathcal{W}_\infty(\Sigma, \gamma)$, but no $\mathcal{C}_0^1$ such function (in fact, not even merely a $\mathcal{C}_0^1$ proper, nor a $\mathcal{C}_0^1$ positive definite, function in $\mathcal{W}(\Sigma, \gamma)$). Specifically, take the following system $\Sigma_1$:

$$\begin{aligned}
\dot{x}_1 &= |x_1|(-x_1 + |x_2| + u_1) \\
\dot{x}_2 &= x_2(-x_1 - |x_2| + u_2) .
\end{aligned}$$

Note that this system is of the form (1), with $n = m = 2$. Now consider the following function (basically, the $L^1$ norm):

$$V_1(x) = 2|x_1| + 2|x_2| .$$

This function is proper and positive definite, and globally Lipschitz. We prove in Section 3:

**Proposition 2.2** For the above system $\Sigma_1$, and unit gain $\gamma = 1$, we have $V_1 \in \mathcal{W}_\infty(\Sigma_1, 1)$. On the other hand, if $V$ is any $\mathcal{C}_0^1$ function in $\mathcal{W}(\Sigma_1, 1)$, then $V$ is not positive definite nor proper.

**Remark 2.3** We chose to give example $\Sigma_1$ because of its simplicity. However, at the cost of added complexity, one may easily provide similar examples which are such that the vector fields $g_i$ are not merely Lipschitz. For a $\mathcal{C}^1$ example, we may take

$$\begin{aligned}
\dot{x}_1 &= |x_1 x_2|^3 + x_1(-|x_1| + u_1) \\
\dot{x}_2 &= -(x_1 x_2)^3 + x_2(-|x_2| + u_2) .
\end{aligned}$$

The proof of the analog of Proposition 2.2, using the same function $V_1$ and gain 1, is virtually the same. □

The second result produces an example showing that there may exist a continuous $V \in \mathcal{W}_\infty(\Sigma, \gamma)$, but every locally Lipschitz function in $\mathcal{W}(\Sigma, \gamma)$ must be either nonproper or non-positive definite. Specifically, we will consider the following system $\Sigma_3$:

$$\begin{aligned}
\dot{x}_1 &= -x_1 + x_2 + u_1 \\
\dot{x}_2 &= 3x_2^{4/3}(-x_1 - x_2 + u_2)
\end{aligned}$$

which is of the form (1), with $n = m = 2$. We will consider the following function:

$$V_3(x_1, x_2) := x_1^2 + x_2^{2/3}$$

which is proper, positive definite, and continuous. We prove in Section 3:

**Proposition 2.4** For the above system $\Sigma_3$, and unit gain $\gamma = 1$, we have $V_3 \in \mathcal{W}_\infty(\Sigma_3, 1)$. On the other hand, if $V$ is any locally Lipschitz function in $\mathcal{W}(\Sigma_3, 1)$, then $V$ is not positive definite nor proper.

There is a general positive result as well. We show that, for any system $\Sigma$ as in (1), and any $\gamma > 0$, $\mathcal{W}_\infty(\Sigma, \gamma) \neq \emptyset$ implies that $\mathcal{W}_\infty(\Sigma, \gamma') \bigcap \mathcal{C}_0^1 \neq \emptyset$ for each $\gamma' > \gamma$. In other words, it is always possible to smoothly approximate a proper positive definite continuous $V$ by one that is continuously differentiable away from zero (actually, the proof provides an infinitely differentiable such approximation), provided that we allow a negligible increase in gain. This is summarized in the following statement, which we also prove in Section 3:

**Theorem 1** *For any system $\Sigma$ as in (1),*

$$\inf \{\gamma \mid \mathcal{W}_\infty(\Sigma, \gamma) \neq \emptyset\} \;=\; \inf \left\{\gamma \mid \mathcal{W}_\infty(\Sigma, \gamma) \bigcap \mathcal{C}_0^1 \neq \emptyset\right\}.$$

This result is significant in so far as the "inf" in question is the one of interest in $H_\infty$ control problems.

**Remark 2.5** The result in Theorem 1 is stated for input-affine systems. It is false in general for arbitrary systems (2). One elegant statement can be made by considering systems of the following special form, for any $p \geq 1$:

$$\dot{x} \;=\; g_0(x) + \sum_{i=1}^m u_i^p \, g_i(x) \tag{8}$$

(again assuming that the vector fields $g_i$, $i = 0, \ldots, m$ are locally Lipschitz in $x$). We may interpret the powers $u^p$, for negative $u$, either as "$|u|^p$" or as "sign $u \, |u|^p$"; with either interpretation, we shall show, also in Section 3, that Theorem 1 holds (for any system of this form) if $p \leq 2$ and does not hold (for some systems of this form) if $p > 2$. □

Finally, we will analyze in Section 4 the special case $m = n = 1$. We will show that there is in that case no gap between the locally Lipschitz and the differentiable case, for systems affine in inputs, but that a gap reappears if we deal with systems that are not input-affine.

## 3 Proofs of Main Results

We first prove Proposition 2.2 and Proposition 2.4. We then prove Theorem 1 (in somewhat more generality), and, finally, we justify the claims made in Remark 2.5.

### 3.1 Proof of Proposition 2.2

In order to verify that $V_1 \in \mathcal{W}(\Sigma_1, 1)$, we compute its subgradients. At any point $x = (x_1, x_2) \in \mathbb{R}^2$ with $x_1 x_2 \neq 0$, obviously

$$\partial_D V_1(x) \;=\; \{\nabla V_1(x)\} \;=\; \left\{2\left(\frac{x_1}{|x_1|}, \frac{x_2}{|x_2|}\right)\right\},$$

so (7) becomes

$$\nabla V_1(x) \cdot \left(g_0(x) + \sum_{i=1}^m u_i g_i(x)\right) = 2\frac{x_1}{|x_1|}|x_1|(-x_1 + |x_2| + u_1) + 2\frac{x_2}{|x_2|}x_2(-x_1 - |x_2| + u_2)$$
$$= 2u_1 x_1 + 2u_2|x_2| - 2x_1^2 - 2x_2^2$$
$$\leq |u|^2 - |x|^2,$$

as desired. Suppose now that $x_1 = 0$ and $x_2 \neq 0$. Then, it is an easy exercise with the definition of subgradients to see that

$$\partial_D V_1(x) = [-2, 2] \times \left\{2\frac{x_2}{|x_2|}\right\},$$

so in this case, for any $\zeta = (\zeta_1, \zeta_2) \in \partial_D V_1(x)$, we have:

$$\zeta \cdot \left(g_0(x) + \sum_{i=1}^m u_i g_i(x)\right) = \zeta_1 \cdot 0 + 2\frac{x_2}{|x_2|}x_2(-|x_2| + u_2)$$
$$= 2u_2|x_2| - 2x_2^2$$
$$\leq |u_2|^2 - |x_2|^2 \leq |u|^2 - |x|^2,$$

again as desired. (The actual form of $\zeta_1$ turns out to be irrelevant.) The case $(x_1 \neq 0, x_2 = 0)$ is similar. So $V_1 \in \mathcal{W}(\Sigma_1, 1)$.

We next show that, if $W \in \mathcal{W}(\Sigma_1, 1)$ is of class $\mathcal{C}_0^1$, then it cannot be proper nor positive definite. So, assume that

$$\nabla W(x) \cdot \left(g_0(x) + \sum_{i=1}^m u_i g_i(x)\right) \leq |u|^2 - |x|^2 \quad (9)$$

holds for all $x \neq 0$ and all $u$. Fix any number $a > 0$. Consider any point of the form $x = (x_1, x_2) = (a, x_2)$, with $x_2 \neq 0$. With the special choice $u_1 = x_1$ and $u_2 = |x_2|$, we have $|x| = |u|$, so (using subscripts to denote partial derivatives) inequality (9) reduces to:

$$W_{x_1}(a, x_2) - (\text{sign } x_2) W_{x_2}(a, x_2) \leq 0.$$

When $x_2 > 0$ this gives $W_{x_1}(a, x_2) - W_{x_2}(a, x_2) \leq 0$, so taking the limit as $x_2 \to 0^+$, we conclude

$$W_{x_1}(a, 0) - W_{x_2}(a, 0) \leq 0$$

for all $a > 0$. Arguing with negative $x_2$ and taking $x_2 \to 0^-$, we get also

$$W_{x_1}(a, 0) + W_{x_2}(a, 0) \leq 0$$

for all $a > 0$. We conclude that $W_{x_1}(a, 0) \leq 0$ for all $a > 0$. This means that $W(\cdot, 0)$ must be bounded above (and is nonnegative on $\mathbb{R}_{\geq 0}$), so $W$ cannot be proper. On the other hand, if $W(0) = 0$, this implies $W(a, 0) \equiv 0$, so $W$ cannot be positive definite either. This completes the proof of Proposition 2.2.

## 3.2 Proof of Proposition 2.4

The intuitive idea of the construction of $\Sigma_2$ and $V_2$ is as follows. We start with $\widetilde{V}_2(\tilde{x}_1, \tilde{x}_2) = \tilde{x}_1^2 + \tilde{x}_2^2$, the square norm in $\mathbb{R}^2$, $(\tilde{x}_1, \tilde{x}_2)$ denoting the canonical coordinates on $\mathbb{R}^2$, and consider the harmonic oscillator motion, but with rescaled time so that the $\tilde{x}_1$-axis consists of equilibria: $\tilde{g}(\tilde{x}_1, \tilde{x}_2) = (\tilde{x}_2 \tilde{x}_2^2, -\tilde{x}_1 \tilde{x}_2^2)$. Note that $\widetilde{V}_2$ is an integral for the motions of $\tilde{g}$. Next, we make a change of coordinates which is a global homeomorphism but fails to be a diffeomorphism: $\tilde{x}_1 \mapsto \tilde{x}_1 =: x_1, \tilde{x}_2 \mapsto \tilde{x}_2^3 =: x_2$. In the new coordinates $(x_1, x_2)$, $V_2$, the transformation of $\widetilde{V}_2$, is no more locally Lipschitz. However, the transformation $g$ of $\tilde{g}$ is locally Lipschitz, and it is of course still true that $V_2$ integrates $g$. The key fact is this: for any locally Lipschitz $W$ which is nonincreasing along trajectories of $g$, $W$ cannot be positive definite nor proper. (In other words, there cannot be any locally Lipschitz "weak Lyapunov function" for $\dot{x} = g(x)$). Finally, we add an input-dependent term to $g$ which provides the dissipation property. When $|u| = |x|$, this dissipation property reduces to the nonincreasing property mentioned above, and hence leads to a contradiction for $V$ locally Lipschitz.

We define the following locally Lipschitz vector field $g$ in $\mathbb{R}^2$:

$$g(x_1, x_2) = \begin{pmatrix} x_2 \\ -3x_1 x_2^{4/3} \end{pmatrix}$$

and note that:

$$\zeta \cdot g(x) = 0 \quad \forall \zeta \in \partial_D V_2(x), \ \forall x \in \mathbb{R}^2. \tag{10}$$

Indeed, pick any $x = (x_1, x_2)$. If $x_2 \neq 0$ then necessarily $\zeta = (2x_1, (2/3)x_2^{-1/3})$, so the claim is clear. If instead $x_2 = 0$, then $g(x_1, x_2) = 0$, so the claim holds as well.

Note that system $\Sigma_2$ is of the form $\dot{x} = f(x, u)$, where $f(x, u) = g(x) + h(x, u)$, with

$$h(x, u) := \begin{pmatrix} u_1 - x_1 \\ 3x_2^{4/3}(u_2 - x_2) \end{pmatrix}.$$

Now note that:

$$\zeta \cdot h(x, u) \leq |u|^2 - |x|^2 \quad \forall \zeta \in \partial_D V_2(x), \ \forall x \in \mathbb{R}^2, \forall u \in \mathbb{R}^2. \tag{11}$$

To verify this, pick any $x = (x_1, x_2)$. Take first the case $x_2 \neq 0$. Since $\zeta$ must be equal to $(2x_1, (2/3)x_2^{-1/3})$, we have that

$$\zeta \cdot h(x, u) = 2x_1 u_1 - 2x_1^2 + 2x_2 u_2 - 2x_2^2 \leq |u|^2 - |x|^2.$$

If instead $x_2 = 0$, then $\zeta = (2x_1, \zeta_2)$, for some $\zeta_2$, so

$$\zeta \cdot h(x, u) = 2x_1 u_1 - 2x_1^2 + \zeta_2 \cdot 0 \leq u_1^2 - x_1^2 \leq |u|^2 - |x|^2.$$

Together with (10), we have thus that $\zeta \cdot f(x, u) \leq |u|^2 - |x|^2$, so we have proved that $V_2 \in \mathcal{W}_\infty(\Sigma_2, 1)$. Now we show that we cannot have a locally Lipschitz $V \in \mathcal{W}_\infty(\Sigma_2, 1)$.

For any $V \in \mathcal{W}(\Sigma_2, 1)$, setting $u = x$ in the inequality $\zeta \cdot f(x, u) \leq |u|^2 - |x|^2$, and using that $h(x, x) = 0$, gives that

$$\zeta \cdot g(x) \leq 0 \quad \forall \zeta \in \partial_D V(x), \ \forall x \in \mathbb{R}^2. \tag{12}$$

We will show that no locally Lipschitz function $V$ can satisfy such a property and be positive definite or proper. Suppose given such a $V$; we will prove that, for any positive number $a$ at which the locally Lipschitz function $x_1 \mapsto V(x_1, 0)$ is differentiable, necessarily $V_{x_1}(a, 0) \leq 0$. This means that $V$ decreases along the $x_1$-axis, and the negative conclusion follows. Fix now any such $a$, and denote $\xi := (a, 0)'$. We will analyze the behavior of $V$ in a neighborhood of $\xi$, by considering a curve which approaches $\xi$ along an orbit of $g$.

Take the following parameterized curve:

$$\gamma : [0, 1] \to \mathbb{R}^2 : t \mapsto \begin{pmatrix} at \\ (a^2 - (at)^2)^{3/2} \end{pmatrix}$$

which has the property that $\gamma(1) = \gamma'(1) = \xi$. Note also that, for each $t \in [0, 1]$,

$$g(\gamma(t)) = \begin{pmatrix} (a^2 - (at)^2)^{3/2} \\ -3at(a^2 - (at)^2)^2 \end{pmatrix} = \beta(t)\, \gamma'(t)$$

where $\beta(t) = \frac{1}{a}(a^2 - (at)^2)^{3/2} > 0$ for all $t < 1$.

**Lemma 3.1** The function $t \in [0, 1] \mapsto W := V(\gamma(t))$ is nonincreasing.

*Proof.* Since $W$ is locally Lipschitz, its derivative exists almost everywhere. We must prove that $W'(t) \leq 0$ for almost all $t$. This will follow from the following statement, valid for all subgradients:

$$\eta \leq 0 \quad \text{for all } t_0 \in [0, 1) \text{ and all } \eta \in \partial_D W(t_0). \tag{13}$$

The idea of the proof is as follows. If $\nabla V(\gamma(t_0))$ exists, then

$$W'(t_0) = \nabla V(\gamma(t_0))\, \gamma'(t_0) = \frac{1}{\beta(t_0)} \nabla V(\gamma(t_0))\, g(\gamma(t_0)) \leq 0$$

where the last inequality follows from (12). However, there is no reason for $V$ to be differentiable at the points in the image of $\gamma$. So we need to apply the "approximate chain rule" for subgradients. Pick any $t_0 \in [0, 1)$ and $\eta \in \partial_D W(t_0)$. Take any $\varepsilon > 0$, and let $0 < \delta < \varepsilon$ be such that

$$|\zeta|\, |g(x) - g(\gamma(t_0))| \leq \varepsilon \beta(t_0) \tag{14}$$

for all $\zeta \in \partial_D V(x)$ and all $x$ such that $|x - \gamma(t_0)| < \delta$ (using continuity of $g$, and noting that all such $\zeta$ are bounded by a Lipschitz constant for $V$ on a ball of radius $\delta$ about $\gamma(t_0)$) and also so that

$$|\zeta|\, |\gamma'(t_1) - \gamma'(t_0)| \leq \varepsilon \tag{15}$$

whenever $\zeta \in \partial_D V(x)$ for any $x$ as above and $|t_1 - t_0| < \delta$ (using now the fact that $\gamma \in \mathcal{C}^1$).

We apply to $\eta$ and $\delta$ the chain rule in [6], Theorem 2.2.5, but in its viscosity rather than subgradient form (the viscosity form follows from the version given in that reference by applying the approximation theorem in [6], Proposition 3.4.5). This tells us that there exist $t_1$, $x$, $\zeta$, and $\zeta'$ such that $|t_1 - t_0| < \delta$, $|x - \gamma(t_0)| < \delta$, $\zeta \in \partial_D V(x)$, and $\zeta' \in \partial_D(\zeta \cdot \gamma)(t_1)$, so that $|\gamma(t_1) - \gamma(t_0)| < \varepsilon$ and $|\eta - \zeta'| < \varepsilon$.

Since $\gamma$ is differentiable, the scalar map $\zeta \cdot \gamma$ is too, and hence $\zeta' = (d/dt)(\zeta \cdot \gamma)(t_1) = \zeta \cdot \gamma'(t_1)$, so

$$\begin{aligned}
\zeta' &= \zeta \cdot \gamma'(t_1) = \zeta \cdot (\gamma'(t_1) - \gamma'(t_0)) + \zeta \cdot \gamma'(t_0) \\
&= \zeta \cdot (\gamma'(t_1) - \gamma'(t_0)) + \frac{1}{\beta(t_0)} \zeta \cdot g(\gamma(t_0)) \\
&= \theta + \frac{1}{\beta(t_0)} \zeta \cdot g(x),
\end{aligned}$$

where
$$\theta = \zeta \cdot (\gamma'(t_1) - \gamma'(t_0)) + \frac{1}{\beta(t_0)} \zeta \cdot (g(\gamma(t_0)) - g(x)) ,$$

and therefore $\zeta' \leq 2\varepsilon$ by (12), (14), and (15). So $\eta \leq 3\varepsilon$. As $\varepsilon$ was arbitrary, we conclude $\eta \leq 0$. ∎

Since $W(t) \geq W(1) = V(\xi)$ for all $t < 1$, for any $h \geq 0$ the last term in the following expression is nonpositive:

$$\frac{1}{h}(V(a,0) - V(a-ah,0)) = \frac{1}{h}(V(\gamma(1-h)) - V(a-ah,0)) + \frac{1}{h}(W(1) - W(1-h)) . \quad (16)$$

Since $V$ is locally Lipschitz, there is some constant $L$ such that (evaluating $\gamma(1-h)$):

$$\frac{1}{h} |V(\gamma(1-h)) - V(a-ah,0)| \leq La^3 h^{\frac{1}{2}} (2-h)^{3/2} \to 0$$

as $h \to 0^+$. We conclude, by taking limits as $h \to 0^+$ in (16), that $V_{x_1}(a,0) \leq 0$, as claimed. ∎

### 3.3 Proof of Theorem 1

In order to be able to justify Remark 2.5, we will prove Theorem 1 for a more general class of systems, namely those of the following general form:

$$\dot{x} = g_0(x) + \sum_{i=1}^m \varphi(u_i) g_i(x) \quad (17)$$

where $\varphi(r) = |r|^p$ or $\varphi(r) = (\operatorname{sign} r) |r|^p$ (interpreting $\varphi(0) = 0$), and $p$ is a fixed real number in the closed interval $[1, 2]$. Note that Theorem 1 as stated would correspond to $p = 1$ and the second choice of $\varphi$. We still suppose that states $x$ evolve in $\mathbb{R}^n$, and inputs $u = (u_1, \ldots, u_m)$ take values in $\mathbb{R}^m$, and that the vector fields $g_i$, $i = 0, \ldots, m$ are locally Lipschitz in $x$.

The proof will be based upon the following general technical fact:

**Lemma 3.2** Assume that we are given:

- an open subset $\mathcal{O}$ of $\mathbb{R}^n$;
- a continuous function $\alpha : \mathcal{O} \to \mathbb{R}_{>0}$;
- a continuous function $\beta : \mathcal{O} \to \mathbb{R}_{\geq 0}$;
- an $\varepsilon > 0$;
- a continuous function $V : \mathcal{O} \to \mathbb{R}_{>0}$ satisfying

$$\zeta \cdot \left(g_0(x) + \sum_{i=1}^m \varphi(u_i) g_i(x)\right) \leq -\alpha(x) + \beta(x) |u|^2 \quad \forall \zeta \in \partial_D V(x) \quad (18)$$

for all $x \in \mathcal{O}$ and $u \in \mathcal{U}$.

Then, there exists a smooth $W : \mathcal{O} \to \mathbb{R}$ such that

$$|V(x) - W(x)| \leq \frac{1}{2} V(x) \tag{19}$$

for all $x \in \mathcal{O}$, and

$$\nabla W(x) \cdot \left( g_0(x) + \sum_{i=1}^{m} \varphi(u_i) \, g_i(x) \right) \leq -\alpha(x) + [(1+\varepsilon)\beta(x) + \varepsilon] \, |u|^2 \tag{20}$$

for all $x \in \mathcal{O}$ and $u \in \mathcal{U}$.

Before proving the lemma, we explain how to obtain Theorem 1 as a corollary. We pick a system $\Sigma$ as in (1), and any $\gamma' > \gamma > 0$. We suppose given a continuous proper and positive definite $V \in \mathcal{W}(\Sigma, \gamma)$, and need to show the existence of some $W \in \mathcal{W}(\Sigma, \gamma')$ which is $\mathcal{C}^1$ on $\mathcal{O} = \mathbb{R}^n \setminus \{0\}$, in addition to being proper and positive definite. We pick, in Lemma 3.2, $\beta(x) \equiv \gamma$, $\alpha(x) = |x|^2$, and any $\varepsilon > 0$ such that $(1+\varepsilon)\gamma + \varepsilon < \gamma'$. The Lemma then applies to the restriction of $V$ to $\mathcal{O}$. We obtain a $W$ as in the Lemma, and extend it from $\mathcal{O}$ to $\mathbb{R}^n$ by defining $W(0) := 0$. The approximation property (19) insures that $W$ is proper and positive definite, and it is continuous at zero as well.

The proof of Lemma 3.2 will be based upon a reduction to the following result, which we quote from [15]:

**Lemma 3.3** Let $\Sigma : \dot{x} = f(x, d)$ be a system, with $x \in \mathbb{X} = \mathbb{R}^n$, and $d \in \mathcal{U}$, a compact metric space, so that $f(x, d)$ is locally Lipschitz in $x$ uniformly on $d$ and jointly continuous in $x$ and $d$. Assume that we are given:

- an open subset $\mathcal{O}$ of $\mathbb{X}$;

- a continuous, nonnegative function $V : \mathcal{O} \to \mathbb{R}$ satisfying

$$\zeta \cdot f(x, d) \leq \Theta(x, d) \quad \forall x \in \mathcal{O}, \, \zeta \in \partial_D V(x), \, d \in \mathcal{U} \tag{21}$$

  with some continuous function $\Theta : \mathcal{O} \times \mathcal{U} \to \mathbb{R}$;

- two positive, continuous functions $\Upsilon_1$ and $\Upsilon_2$ on $\mathcal{O}$.

Then, there exists a smooth $\hat{V} : \mathcal{O} \to \mathbb{R}$ such that

$$\left| V(x) - \hat{V}(x) \right| \leq \Upsilon_1(x) \quad \forall x \in \mathcal{O} \tag{22}$$

and

$$\nabla \hat{V}(x) \cdot f(x, d) \leq \Theta(x, d) + \Upsilon_2(x) \quad \forall x \in \mathcal{O}, \, d \in \mathcal{U}. \tag{23}$$

$\square$

The proof given in [15] employs tools from nonsmooth analysis, borrowing in particular several simple facts from [5, 6, 18] in order to pass from continuous $V$ to locally Lipschitz $V$, followed by a standard smoothing argument as given in [16]. (To be precise, the result is stated and proved in [15] using the language of proximal subgradients rather than viscosity subgradients. The proof, however, is the same in both cases.)

**Proof of Lemma 3.2**

We define
$$f(x,d) := \left(1 - |d|^2\right) g_0(x) + \sum_{i=1}^m \varphi(d_i) \left(1 - |d|^2\right)^{1-\frac{p}{2}} g_i(x)$$

with $\mathcal{U}$ = unit ball in $\mathbb{R}^m$. (If $p = 2$, we interpret the coefficient of $g_i$ as simply $\varphi(d_i)$.) Note that this system satisfies the hypotheses of Lemma 3.3. We use the same $\mathcal{O}$ and $V$, and let
$$\Theta(x,d) := -\left(1 - |d|^2\right)\alpha(x) + |d|^2 \beta(x).$$

We need to verify that (21) indeed holds. Pick any $x$ and $\zeta$, and any $d$ with $|d| \leq 1$. We treat separately the cases $|d| = 1$ and $|d| < 1$.

*Case 1:* $|d| = 1$. If $p < 2$, then (21) holds trivially (because $f(x,d) = 0$ and $\Theta(x,d) = \beta(x) \geq 0$). So suppose $p = 2$. We must verify that
$$\sum_{i=1}^m \varphi(d_i)\, \zeta \cdot g_i(x) \leq \beta(x) \tag{24}$$

where $\varphi(r) = r^2$, or $\varphi(r) = r|r|$, and $d$ is a vector of unit norm. We are assuming that (18) holds, that is,
$$\zeta \cdot g_0(x) + \sum_{i=1}^m \varphi(u_i)\, \zeta \cdot g_i(x) \leq -\alpha(x) + \beta(x)|u|^2 \tag{25}$$

for all $u$. Now pick any $u$ of the form $(1/\varepsilon)d$, with $\varepsilon > 0$, so $\varepsilon^2 \varphi(u_i) = \varphi(d_i)$ for each $i$. Multiplying both sides of (25) by $\varepsilon^2$, we have that
$$\varepsilon^2 \zeta \cdot g_0(x) + \sum_{i=1}^m \varphi(d_i)\, \zeta \cdot g_i(x) \leq -\varepsilon^2 \alpha(x) + \beta(x),$$

so taking limits as $\varepsilon \to 0$ one obtains the desired inequality (24).

*Case 2:* $|d| < 1$. Introduce, for such a $d$, $u = (u_1, \ldots, u_m)$, where:
$$u_i := \frac{d_i}{\sqrt{1 - |d|^2}} \quad i = 1, \ldots, m.$$

Observe that
$$|u|^2 = \frac{|d|^2}{1 - |d|^2}, \quad 1 - |d|^2 = \frac{1}{1 + |u|^2}, \quad \varphi(d_i)\left(1 - |d|^2\right)^{1-\frac{p}{2}} = \frac{\varphi(u_i)}{1 + |u|^2} \tag{26}$$

for all $i$, so that
$$f(x,d) = \frac{1}{1 + |u|^2}\left(g_0(x) + \sum_{i=1}^m \varphi(u_i)\, g_i(x)\right)$$

and thus (18) implies that $\zeta \cdot f(x,d) \leq \Theta(x,d)$ for all $\zeta \in \partial_D V(x)$, as wanted.

So, we apply Lemma 3.3 with
$$\Upsilon_1(x) := \frac{1 - \delta}{4} V(x)$$

and
$$\Upsilon_2(x) := \delta \min\{1, \alpha(x)\}$$
where $\delta \in (0,1)$ is picked such that
$$\frac{1}{1-\delta} \leq 1+\varepsilon \quad \text{and} \quad \frac{\delta}{1-\delta} \leq \min\left\{\varepsilon, \frac{1}{4}\right\}.$$

This provides a function $\hat{V}$.

Now fix any $u \in \mathbb{R}^m$, and introduce
$$d_i := \frac{u_i}{\sqrt{1+|u|^2}} \quad i = 1, \ldots, m.$$

Observe that the relations (26) hold, and hence also
$$|d|^2 = \frac{|u|^2}{1+|u|^2} < 1$$
holds. Furthermore, $g_0(x) + \sum_{i=1}^{m} \varphi(u_i) g_i(x) = (1+|u|^2) f(x,d)$. Therefore, (23) gives

$$\begin{aligned}
\nabla \hat{V}(x) \cdot \left( g_0(x) + \sum_{i=1}^{m} \varphi(u_i) g_i(x) \right) \\
\leq \ & (1+|u|^2) \left( -\left(1-|d|^2\right) \alpha(x) + |d|^2 \beta(x) + \delta \min\{1, \alpha(x)\} \right) \\
= \ & -\alpha(x) + |u|^2 \beta(x) + \delta \min\{1, \alpha(x)\} + |u|^2 \delta \min\{1, \alpha(x)\} \\
\leq \ & -\alpha(x) + |u|^2 \beta(x) + \delta \alpha(x) + |u|^2 \delta \\
= \ & -(1-\delta)\alpha(x) + |u|^2 (\beta(x) + \delta)
\end{aligned}$$

and if we now let
$$W := \frac{1}{1-\delta} \hat{V}$$
we conclude that
$$|V - W| \leq \frac{\delta}{1-\delta} V + \frac{1}{1-\delta} \left| V - \hat{V} \right| \leq \frac{1}{2} V,$$
and
$$\nabla W(x) \cdot \left( g_0(x) + \sum_{i=1}^{m} \varphi(u_i) g_i(x) \right) \leq -\alpha(x) + |u|^2 \frac{\beta(x) + \delta}{1-\delta} \leq -\alpha(x) + [(1+\varepsilon)\beta(x) + \varepsilon] |u|^2$$

as desired for the conclusions of Lemma 3.2. ∎

Notice, incidentally, that $W$ is smooth ($\mathcal{C}^\infty$), not merely $\mathcal{C}^1$.

### 3.4 Results for systems of form (8)

The positive claim in Remark 2.5, for $p \leq 2$, has been proved already. To show that Theorem 1 does not extend to systems (17) when $p > 2$, it will be convenient to introduce the vector field corresponding to what one might call an "$L^1$ harmonic oscillator". This is the locally Lipschitz vector field on $\mathbb{R}^2$ defined by

$$g(x_1, x_2) := \begin{pmatrix} |x_1| x_2 \\ -|x_2| x_1 \end{pmatrix}. \tag{27}$$

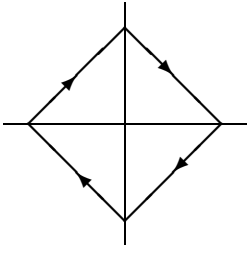

Figure 1: The vector field $g$

The flow of $g$ leaves invariant the $L^1$-balls around the origin, see Figure 1. Associated to $g$ is the $L^1$ norm, seen as a proper, positive definite Lipschitz function

$$V_1(x_1, x_2) := |x_1| + |x_2|.$$

Except for the fact that we do not now need the factor "2", this function is as before. So $\partial_D V_1(x) = \{(x_1/|x_1|, x_2/|x_2|)\}$ when $x_1 x_2 \neq 0$, $\partial_D V_1(x) = [-1, 1] \times \{x_2/|x_2|\}$ when $x_1 = 0$ and $x_2 \neq 0$, and similarly if $x_2 = 0$. Therefore,

$$\zeta \cdot g = 0 \quad \text{for all} \ \zeta \in \partial_D V_1(x), \ x \in \mathbb{R}^n. \tag{28}$$

We also introduce the following vector field:

$$g_0(x_1, x_2) := \begin{pmatrix} -|x_1| x_1 \\ -|x_2| x_2 \end{pmatrix}. \tag{29}$$

For the same $V_1$,

$$\zeta \cdot g_0 = -|x|^2 \quad \text{for all} \ \zeta \in \partial_D V_1(x), \ x \in \mathbb{R}^n. \tag{30}$$

Now fix any $p > 2$, and consider the following system $\Sigma_p$, which is defined in terms of the above vector fields:

$$\dot{x} = g_0(x) + |u_1|^p g(x) - |u_2|^p g(x) \tag{31}$$

This is a system of type (8) with $m = 2$, $g_1 = g$, and $g_2 = -g$. In view of (28) and (30), we have that

$$\zeta \cdot (g_0(x) + |u_1|^p g(x) - |u_2|^p g(x)) \leq -|x|^2$$

for all $x \in \mathbb{R}^2$ and $u = (u_1, u_2) \in \mathbb{R}^2$, and every $\zeta \in \partial_D V_1(x)$. This means that $V_1 \in \mathcal{W}_\infty(\Sigma_p, \gamma)$, for any $\gamma > 0$. We now see that the equality in Theorem 1 cannot hold. In fact, we prove the following far stronger statement:

**Proposition 3.4** For all $V \in \mathcal{C}_0^1$, and all $\gamma > 0$, $V \notin \mathcal{W}(\Sigma_p, \gamma)$.

*Proof.* Suppose that $V \in \mathcal{C}_0^1$ and

$$\nabla V(x) \cdot (g_0(x) + |u_1|^p g(x) - |u_2|^p g(x)) \leq -|x|^2 + \gamma |u_1|^2 + \gamma |u_2|^2 \tag{32}$$

for all $x \neq 0$ and all $u = (u_1, u_2)$. Fix any $x = \xi \neq 0$. Taking $u = 0$ gives that $\nabla V(\xi) \neq 0$. Since $p > 2$, letting $u_2 = 0$ and $u_1 \to +\infty$ gives $\nabla V(\xi) \cdot g(\xi) \leq 0$. On the other hand, $u_1 = 0$ and $u_2 \to +\infty$ gives $-\nabla V(\xi) \cdot g(\xi) \leq 0$. So we conclude that $\nabla V(\xi) \cdot g(\xi) = 0$ for all $\xi \neq 0$.

Now fix an $a > 0$ and suppose that $\xi$ has the form $(a, x_2)$ with $x_2 \neq 0$. Then, $\nabla V(\xi) \cdot g(\xi) = 0$ means that
$$(\text{sign } x_2) V_{x_1}(a, x_2) - V_{x_2}(a, x_2) = 0$$
Letting separately $x_2 \to 0^+$ and $x_2 \to 0^-$ gives $V_{x_1}(a, 0) = V_{x_2}(a, 0) = 0$, contradicting $\nabla V(a, 0) \neq 0$. ∎

A variation of this example is as follows. Fix again any $p > 2$. Instead of (31), we consider now the following system $\Sigma'_p$:
$$\dot{x} = g_0(x) + \text{sign } u \, |u|^p \, g(x) \tag{33}$$
with $m = 1$. It is also true that $V_1 \in \mathcal{W}_\infty(\Sigma_p, \gamma)$, for any $\gamma > 0$, and Proposition 3.4 again holds, simply taking the separate limits as $u \to +\infty$ or $u \to -\infty$ in its proof.

**Remark 3.5** There is nothing very special about the form (8). Mainly, we picked these systems in order to illustrate with a specific class when the theorem holds. But Theorem 1 holds also for a general class of systems of subquadratic growth in inputs. Specifically, we may consider systems of the general form $\dot{x} = g(x, u)$, where $g$ is jointly continuous in $x$ and $u$, and satisfies, for some constant $p \in (0, 2)$ and some continuous function $h : \mathbb{R}^n \to \mathbb{R}$:

1. $|g(x, u)| \leq h(x) \, (1 + |u|^p)$ for every $x$ and $u$;

2. $|g(x, u) - g(y, u)| \leq c|x - y| \, (1 + |u|^p)$ for every $u$, and every $x, y$ in some (arbitrary chosen) ball.

The proof is the same as for systems of the form (8). One only needs to set
$$f(x, d) = \left(1 - |d|^2\right) g\left(x, d\left(1 - |d|^2\right)^{-1/2}\right)$$
provided $|d| < 1$, and $f(x, d) = 0$ otherwise, and to replace everywhere $g_0(x) + \sum_{i=1}^m \varphi(u_i) g_i(x)$ by $g(x, u)$. We omit the simple details. □

## 4 One-Dimensional Systems

For one-dimensional systems, there is no gap between the locally Lipschitz and the $\mathcal{C}_0^1$ case:

**Proposition 4.1** Suppose given a system $\Sigma$ as in (1), with $n = 1$. Assume that for some $\gamma > 0$ there exists a locally Lipschitz $V \in \mathcal{W}_\infty(\Sigma, \gamma)$. Then $\mathcal{W}_\infty(\Sigma, \gamma) \bigcap \mathcal{C}_0^1 \neq \emptyset$.

*Proof.* Pick a $V$ as in the statement of the proposition. By Rademacher's theorem there exist a zero measure (Borelian) set $N \subset \mathbb{R}$ and a continuous positive function $h$ such that
$$x \notin N \Rightarrow V'(x) \text{ exists and } |V'(x)| \leq \frac{1}{2} h(x). \tag{34}$$

We shall construct $W$ only on $\mathbb{R}_{\geq 0}$, the construction being similar on $\mathbb{R}_{\leq 0}$. If $x \notin N$, $\partial_D V(x) = \{V'(x)\}$, hence

$$V'(x)\left(g_0(x) + \sum_{i=1}^m u_i g_i(x)\right) \leq \gamma |u|^2 - x^2 \quad \forall u \in \mathbb{R}. \tag{35}$$

Since the system $\dot{x} = g_0(x)$ is (globally) asymptotically stable, $g_0(x) < 0$ for each $x > 0$. It follows from (34) and (35) (with $u = 0$) that

$$\frac{1}{2}h(x) \geq V'(x) \geq \frac{x^2}{|g_0(x)|} > 0 \quad \text{for } x > 0,\ x \notin N. \tag{36}$$

Set for any $x > 0$

$$F(x) := \left\{p \geq 0 :\ p\left(g_0(x) + \sum_{i=1}^m u_i g_i(x)\right) \leq \gamma |u|^2 - x^2 \quad \forall u \in \mathbb{R}\right\}.$$

We claim that the (closed convex) set $F(x)$ is *nonempty* for any $x > 0$. Indeed, if $x \notin N$ $V'(x) \in F(x)$. If $x \in N$ we may pick a sequence $(x_n)$ in $\mathbb{R}_{>0} \setminus N$ such that $x_n \to x$ as $n \to \infty$. Since $|V'(x_n)| \leq \frac{1}{2}h(x_n) \leq c$ for some constant $c$, we may extract a subsequence of $(V'(x_n))$ which converges towards some $p \geq 0$. Clearly $p \in F(x)$.

Let us set

$$a := \sum_{i=1}^m g_i(x)^2,\ b := 4\gamma g_0(x),\ c := 4\gamma x^2$$

and

$$\Delta(p) := ap^2 + bp + c.$$

It is a straightforward exercise to see that for any $p \geq 0$,

$$p \in F(x) \iff \Delta(p) \leq 0. \tag{37}$$

By the claim, $\Delta(p) \leq 0$ for some $p$, whereas $\Delta(p) \to \infty$ as $p \to -\infty$ (since $b < 0$). It follows that $\Delta$ has (at least) one *real* root, hence $b^2 - 4ac \geq 0$. Set

$$p(x) := \begin{cases} h(x) & \text{if } a = 0, \\ \min\left\{h(x), \frac{-b + \sqrt{b^2 - 4ac}}{2a}\right\} & \text{if } a \neq 0. \end{cases} \tag{38}$$

Note that $-b = 4\gamma |g_0(x)| > 0$ on $x > 0$, so when $a \to 0$ the second expression in the minimum above becomes unbounded; from here, it follows that $p$ is continuous on $\mathbb{R}_{>0}$. We claim that

$$V'(x) \leq p(x) \quad \forall x \in \mathbb{R}_{>0} \setminus N. \tag{39}$$

Pick any $x \in \mathbb{R}_{>0} \setminus N$. If $a = 0$, then (39) follows from (36) and (38). If $a \neq 0$, since $V'(x) \in F(x)$, we see that $\Delta(V'(x)) \leq 0$ (by (37)), hence $V'(x) \leq \frac{-b + \sqrt{b^2 - 4ac}}{2a}$ which, combined with (36) and (38), yields $V'(x) \leq p(x)$.

We are now ready to define $W$ on $\mathbb{R}_{\geq 0}$: we set $W(x) = \int_0^x p(s)\,ds$ for any $x \geq 0$. Since $0 \leq p \leq h$, $W$ is a (well-defined) locally Lipschitz function on $\mathbb{R}_{\geq 0}$ which is $C^1$ away from 0. Integrating in (39) we get $W(x) \geq V(x) > 0$ for any $x > 0$ and $\lim_{x \to \infty} W(x) = \infty$. Finally we claim that

$$\Delta(p(x)) \leq 0 \quad \forall x \in \mathbb{R}_{>0} \setminus N. \tag{40}$$

Pick any $x \in \mathbb{R}_{>0} \setminus N$. If $a = 0$, then

$$\begin{aligned}\Delta(p(x)) &\leq \Delta(V'(x)) \quad \text{(using } b<0 \text{ and (39))} \\ &\leq 0 \quad \text{(by (37))}.\end{aligned}$$

If $a \neq 0$, we are led to prove that

$$p(x) \in \left[\frac{-b - \sqrt{b^2 - 4ac}}{2a}, \frac{-b + \sqrt{b^2 - 4ac}}{2a}\right].$$

Owing to the definition of $p(x)$, we only have to prove

$$h(x) \geq \frac{-b - \sqrt{b^2 - 4ac}}{2a}.$$

But

$$\frac{-b - \sqrt{b^2 - 4ac}}{2a} = \frac{2c}{|b| + \sqrt{b^2 - 4ac}} \leq \frac{2c}{|b|} = \frac{2x^2}{|g_0(x)|} \leq h(x),$$

by (36). This completes the proof of (40). Note that (40) also holds true for *any* $x > 0$, by continuity of $\Delta \circ p$. Using (37) it means that $W'(x) = p(x) \in F(x)$ for all $x > 0$, i.e.

$$W'(x) \left(g_0(x) + \sum_{i=1}^{m} u_i g_i(x)\right) \leq \gamma |u|^2 - x^2 \quad \forall u \in \mathbb{R}.$$

The proof of the proposition is complete. ∎

We now show that, for systems not affine in inputs, there is again a gap between the Lipschitz and differentiable cases. In Proposition 2.2 we gave an example of a system $\Sigma_1$ and a Lipschitz $V_1 \in \mathcal{W}_\infty(\Sigma_1, 1)$, such that no $V$ in $\mathcal{W}(\Sigma_1, 1) \cap \mathcal{C}_0^1$ can be positive definite or proper. The system was two-dimensional ($n = 2$) and had two-dimensional input ($m = 2$). We now provide a scalar ($n = m = 1$) system $\Sigma_3$ which has the following properties: (a) it admits a Lipschitz $V_3 \in \mathcal{W}_\infty(\Sigma_3, 1)$, but (b) any $W$ in $\mathcal{W}(\Sigma_3, 1)$ must be non-differentiable. Thus, the conclusions are even stronger than for $\Sigma_1$; on the other hand, $\Sigma_3$ is not affine in inputs.

Before giving the form of $\Sigma_3$, we start by considering the following two functions $\varphi_+$ and $\varphi_- : \mathbb{R}^2 \to \mathbb{R}$:

$$\varphi_+(s, t) := \max\{\min\{(s-t)/2, s\}, 0\}$$

and

$$\varphi_-(s, t) := -\max\{\min\{(-s-t)/2, -s\}, 0\}.$$

Both of these are Lipschitz, and $\varphi_+(0, t) = \varphi_-(0, t) = 0$ for all $t$. Thus, the function

$$\varphi(s, t) := (\text{sign } s) \max\left\{\min\left\{\frac{1}{2}(|s| - t), |s|\right\}, 0\right\}$$

(with $\varphi(0, t) \equiv 0$) obtained by using $\varphi_+$ for $s \geq 0$ and $\varphi_-$ for $s < 0$ is also Lipschitz.

As $\min\{\frac{1}{2}(|s| - t), |s|\} \leq |s|$, it follows that

$$|\varphi(s, t)| \in [0, |s|]$$

for all $(s,t)$. Since $\text{sign}\,\varphi(s,t) = \text{sign}\,s$, this means that $\varphi(s) \in [0, s]$ when $s \geq 0$ and $\varphi(s) \in [s, 0]$ when $s \leq 0$. When $t \geq |s|$, $\min\left\{\frac{1}{2}(|s| - t), |s|\right\} = \frac{1}{2}(|s| - t)$, so

$$\varphi(s, t) = 0 \text{ if } t \geq |s|.$$

If $t \leq -|s|$ then $|s| = |s|/2 + |s|/2 \leq |s|/2 - t/2$, so $\min\left\{\frac{1}{2}(|s| - t), |s|\right\} = |s|$. Therefore

$$\varphi(s, t) = s \text{ if } t \leq -|s|.$$

The diagram in Figure 2 summarizes this information about the range of $\varphi$.

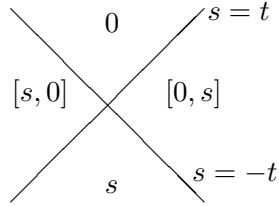

Figure 2: The range of $\varphi$

We define now, for $a, b \geq 0$:

$$\psi(a, b) := \frac{1}{2}\Big[\varphi(b - a, b + a - 2) + (b - a)\Big].$$

As $|b - a| \leq b + a - 2$ if and only if both $a \geq 1$ and $b \geq 1$, $\psi(a, b) = \frac{1}{2}(b - a)$ for such $a, b$. Similarly, the properties of $\varphi$ also imply that $\psi(a, b) = b - a$ when $a \leq 1$ and $b \leq 1$, that $\psi(a, b) \in [b - a, \frac{1}{2}(b - a)]$ when $a \geq b$, and that $\psi(a, b) \in [\frac{1}{2}(b - a), b - a]$ when $a \leq b$. We are now ready to specify $\Sigma_3$. We let:

$$f(x, u) := (|u| + x)\,\psi(x, |u|)$$

for $x \geq 0$, and

$$f(x, u) := x^2 + |u|\,\psi(0, |u|)$$

for $x < 0$. This function is locally Lipschitz. Notice that $\psi(0, |u|) \geq |u|/2 \geq 0$ for all $u$, so that $f(x, u) \geq x^2$ for $x < 0$. Observe the following properties, for all $x \geq 0$:

1. $x \leq 1 \text{ \& } |u| \leq 1 \Rightarrow f(x, u) = u^2 - x^2$.
2. $x \leq 1 \text{ \& } |u| \geq 1 \Rightarrow f(x, u) \leq u^2 - x^2$.
3. $x \geq 1 \text{ \& } |u| \leq 1 \Rightarrow f(x, u) \leq \frac{1}{2}(u^2 - x^2)$.
4. $x \geq 1 \text{ \& } |u| \geq 1 \Rightarrow f(x, u) = \frac{1}{2}(u^2 - x^2)$.

Finally, the Lipschitz function $V_3$ shown in Figure 3.

We show that $V_3 \in \mathcal{W}_\infty(\Sigma_3, 1)$, that is, $V'(x) f(x, u) \leq u^2 - x^2$ for all $x \neq 0$ and all $u$. For $x < 0$ this is obvious, since $V'(x) f(x, u) = -f(x, u) \leq -x^2 \leq u^2 - x^2$ for all $u$. So we only need to analyze the case $x > 0$. For any $0 < x < 1$, $V'(x) f(x, u) = f(x, u) \leq u^2 - x^2$ (properties 1 and 2) while for $x > 1$ we have $V'(x) f(x, u) = 2f(x, u) \leq u^2 - x^2$ as well (properties 3 and 4). Finally, we deal with the nondifferentiability point $x = 1$. An easy calculation shows that $\partial_D V_3(1)$ is the closed interval $[1,2]$. Pick any $\rho \in [1,2]$. When $|u| \leq 1$, $f(x,u) = u^2 - x^2 \leq 0$, so $\rho \geq 1$ implies $\rho f(x,u) = \rho(u^2 - x^2) \leq u^2 - x^2$. Finally, if $|u| \geq 1$, $f(x,u) = \frac{1}{2}(u^2 - x^2) \geq 0$, then $\rho \leq 2$ implies $\rho f(x,u) \leq u^2 - x^2$ as well.

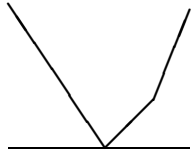

Figure 3: $V_3(x) := \max\{|x|, 2x-1\}$

We now show that any $W$ in $\mathcal{W}(\Sigma_3, 1)$ must be non-differentiable. Suppose that there is some such $W$ which is differentiable. Fix $u = 1$. For any $x$, we must have, since $W \in \mathcal{W}(\Sigma_3, 1)$:

$$W'(x) f(x,1) \leq 1 - x^2.$$

When $x \in (0,1)$, $f(x,1) = 1 - x^2 > 0$, so $W'(x) \leq 1$. On the other hand, for $x > 1$ we have $f(x,1) = (1-x^2)/2 < 0$, so $W'(x) \geq 2$ for such $x$. If $W \in \mathcal{C}_0^1$, this gives a contradiction as $x \to 1^+$ and $x \to 1^-$. However, continuity of $W'$ is not needed for the contradiction, since we can argue as follows. By the mean value theorem,

$$\limsup_{x \to 1^-} \frac{W(x) - W(1)}{x - 1} \leq 1 < 2 \leq \liminf_{x \to 1^+} \frac{W(x) - W(1)}{x - 1},$$

which contradicts the existence of $W'(1)$. ∎

## Acknowledgments

One of the authors (L.R.) thanks the department of mathematics of the Politecnico Di Torino for its hospitality and is deeply grateful to Prof. A. Bacciotti for bringing some nice counterexample in [1] to its attention and for stimulating discussions.